\documentclass{article}[12]
\usepackage{amsmath,amscd,amsthm}
\usepackage{amssymb}
\usepackage{amsfonts}
\usepackage{mathrsfs}
\usepackage{graphics}
\usepackage{epsfig}
\usepackage[top=1.5in,bottom=1in,left=1.25in,right=1.25in]{geometry}

\begin{document}

\title{On John Mather's work}

\author{Sen Hu}
\maketitle

\begin{center}
{\small
School of Mathematics, University of Science and Technology of China\\
Wu Wen-Tsun Key Lab of Mathematics, Chinese Academy of Sciences}
\end{center}

{\tableofcontents}

\section{Mather's work on Singularities}

\subsection{Some backgrounds}

 Singularity theory is to study maps $f: M \rightarrow N$ around singular points up to coordinate changes
 (smooth, holomorphic or topological). Let us start with smooth coordinate changes.

 In \cite{Thom} R. Thom introduced the basic notion of stable maps.

 ${\bf Definition}$: A differentiable mapping $f: M \rightarrow N$ is stable if for any differentiable mapping $\tilde{f}: M \rightarrow N$ sufficiently close to $f$, there are diffeomorphisms $h: M \rightarrow M$ and $k: N \rightarrow N$ such that $h \tilde{f} = f k$.

 He was interested in caustics and the concept of stable maps is the key to understand it.

 Here are some examples of stable maps.

 1) Morse functions,

 2) Whitney's classification of stable maps \cite{Whitney}:
  $f: M^{2} \rightarrow N^{2}$

 regular point: $y_{1} = x_{1}, y_{2} = x_{2}$,

 fold: $y_{1} = x_{1}, y_{2} = x_{2}^{2}$,

 cusp: $y_{1} = x_{1}, y_{2} = x_{1} x_{2} - \frac{1}{3} x_{2}^{3}$.

   According to Whitney \cite{Whitney}: "A fundamental problem is  to determine what sort of singularities any good approximation
   $f$ to $f_{0}$ must have; what sorts of sets they occupy, what $f$ is like near such points, what topological
   properties hold with reference to them, etc." In a series of papers Mather \cite{MS1, MS2, MS3, MS4, MS5, MS6, MS7}
   singly handed established the foundations of singularity theory.

   We will give a brief introduction to Mather's theory. It is strongly recommended to read his original papers which show
   great rigor and clarity in establishing such a theory. Needless to say that Mather's works are on the shoulders of Thom, Whitney
   and Malgrange.

   For example one needs to measure how close two differentiable mappings are. For this Whitney's $C^{k}$ topology is the key.
   Let us first introduce the notion of $k$-jets following Mather.

     The notion of a $k$-jet of a mapping of $N$ into $P$, with source $x$ and target $y$ is an open neighborhood of $y$, the equivalence relation being $k$-tangentiality at $x$.
     
      One denotes the set of all $k$-jets of mappings of $N$ into $P$ with source $x$ and target $y$ by $J^{k}(N, P)_{x, y}$.
       If $z \in J^{k}(N, P)_{x, y}$ is in the equivalence class of $f$, one says $z$ is the $k$-jet of $f$ at $x$ and $f$ is a representative of $z$. One writes $z = j^{k}(f)(x)$. We set $J^{k}(N, P) = \cup_{x \in N, y \in P} J^{k}(N, P)_{x, y}$.

      {\bf The Whitney $C^{k}$ topologies:} Let $N$ and $P$ be manifolds, and $k$ a non-negative integer. For any subset $U$ of $J^{k}(N, P)$,
   where it is the space of $k$-jets of maps from $N$ to $P$,
   
   set

   $$M(U) = \{f \in C^{k}(N, P): j^{k}(f)(N) \subset U\}.$$

    The family of all sets $M(U)$, where $U$ ranges over all subsets of $J^{k}(N, P)$, is the basis of a topology on $C^{k}(N, P).$

    For $C^{\infty}$, we set

    $$W_{\infty} = W \mbox{ denote the topology generated by} \cup_{k} W_{k}$$

    where $k$ runs over all non-negative integers.

    \subsection{Mather-Malgrange's preparation theorem}

 One of the crucial ingredients in singularity theory is B. Malgrange's preparation theorem for $C^{\infty}$ mappings \cite{Malgrange1, Malgrange2}.

 {\bf Malgrange's Preparation Theorem}: Let $f: {\bf R}^{m} \rightarrow {\bf R}^{n}, y = f(x)$ be a smooth map with $0$ its critical point.
 The local ring is $Q(f) = {\bf R}[x_{1}, ..., x_{n}] / \{\frac{\partial f}{\partial x_{1}}, ..., \frac{\partial f}{x_{n}}\}$.

 Suppose that the local ring $Q$ is of finitely dimensional as a real vector space, 
 and let $\{e_{1}(x), ..., e_{r}(x)\}$ in ${\cal E}(x)$ be a basis of $Q$. Then every 
 germ $\phi(x)$ of ${\cal E}(x)$ can be expressed as
 $\phi(x) = \phi_{1}(y)e_{1}(x) + ... + \phi_{r}(y) e_{r}(x)$, that is, 
 ${\cal E}(x)$ is a finite-dimensional ${\cal E}(y)$-module.

It extends Weierstrass preparation theorem from holomporphic category to differentiable category.
Mather \cite{MS1} proved a stronger version of the Malgrange preperation theorem as his division theorem.

  {\bf Mather's Division Theorem}: Let $U$ be open in ${\bf R}^{n}$. Let $f \in C^{\infty}({\bf R} \times U)$, and let
 $u: U \rightarrow {\bf R}^{p}$ be $C^{\infty}$.

 There exist mappings $Q = Q_{p}: C^{\infty} \times  {\bf R} \times {\bf R}^{p} \rightarrow {\bf R}, H = H_{p}: C^{\infty} \times
 {\bf R} \times {\bf R}^{p} \rightarrow {\bf R}^{p}$
 with the following properties:

 a) For all $f \in C^{\infty}({\bf R}), t \in {\bf R}, u \in {\bf R}^{p}$

  $$f(t) = \Gamma_{p}(t, u) Q_{p}(f, t, u) + R_{p}(t, H_{p}(f, u))$$

 b) $Q$ and $H$ are ${\bf R}-$linear in the first variable.

 c) If $f \in C^{\infty}({\bf R}), u \in {\bf R}^{p}$, and $f$ vanishes on $N(u)$, then $H(f, u) = 0$.
 $N(u) = \{t \in {\bf R}, <t, u> \in N \}, N$ any open neighborhood of $Z$ in ${\bf R} \times {\bf R}^{p}$.

 d) For all $f \in C^{\infty}({\bf R})$, the mappings $Q_{f}: {\bf R} \times {\bf R}^{p} \rightarrow {\bf R}, H_{f}: {\bf R} \times
 {\bf R}^{p} \rightarrow {\bf R}$  defined by $Q_{f}(t, u) = Q(f, t, u)$ and $H_{f}(t, u) = H(f, u)$ are $C^{\infty}$.

 e) For any non-negative integer $k$, there exists $K$ such that

  i) For every compact set $X \subset {\bf R}^{p}$, there exists $C > 0$, such that, for every $f \in C^{\infty}({\bf R}),$

  $$ ||H_{f}||^{*}_{k, X} \leq C ||f||^{*}_{K, N(X)}.$$

  ii) Let $\pi: {\bf R} \times {\bf R}^{p} \rightarrow {\bf R}^{p}, \rho: {\bf R} \times {\bf R}^{p} \rightarrow {\bf R}$ be projections.
  For every $Y \subset {\bf R} \times {\bf R}^{p}$ compact, there exist $C > 0$, such that for every $f \in C^{\infty}({\bf R})$,

  $$ ||Q_{f}||^{*}_{k, Y} \leq C ||f||^{*}_{K, Y^{'}},$$

  where $Y^{'} = N(\pi U) \cup \rho Y$.
  
  $||f||^{*}_{k, U}$ is the usual norm for $C^{k}$ maps in the Banach space of maps.

  To prove the theorem he uses Fourier's integral formula in the form

 $$g(t) = \frac{1}{\pi} \int_{0}^{\infty} d \tau \int_{- \infty}^{\infty} dx g(x) \cos \tau (x - t),$$

 which is valid for any $g \in C^{\infty}({\bf R})$ with compact support, and  any $t \in {\bf R}$. From the above representation
 it is enough to solve the division problem for such functions with certain estimates. Lemma 1 in \cite{MS1} provides
 such results and estimates.

 The division theorem implies Malgrange's preparation theorem.

\subsection{Infinitesimal stable implies stable}

 An important consequence of Mather-Malgrange's preparation theorem is that
if a mapping is proper and infinitesimal stable at a point then it is stable at the point.

 A mapping $f: N \rightarrow P$ is said to be proper if $f^{-1}(K)$ is compact for every compact subset $K$ of $P$.

 {\bf Definition: (Infinitesimal stable)} A $C^{\infty}$-map $f:X \rightarrow Y$ is infinitesimally stable if,
 for every $C^{\infty}$-vector field along $f, \xi: X \rightarrow TY$,
 there exist $C^{\infty}$-vector fields $\chi: X \rightarrow TX$ on $X$ and $\eta: Y \rightarrow TY$ on $Y$,
 such that $\xi=(Tf)\circ \chi + \eta \circ f$.

 A vector field along $f$ is a mapping $\omega: N \rightarrow TP$ such that $\omega(n) \in P_{f(n)}$ for all $n \in N$.

 In \cite{MS2} Mather proved that a proper $C^{\infty}$-mapping between $C^{\infty}$-manifolds is stable if and only if it is infinitesimally stable.
 More precisely,

 ${\bf Theorem} \cite{MS2}$: Let $f: N \rightarrow P$ be $C^{\infty}$. Let $M$ be a closed submanifold of $N$ of codimension $0$.
 Suppose $f|M$ is proper and infinitesimally stable. Then there exists a neighborhood $U$ of $f|M$ in
 $C^{\infty}(M, P)$ and continuous mappings $H_{1}: U \rightarrow Diff^{\infty}(N)$ and $H_{2}: U \rightarrow Diff^{\infty}(P)$
 such that $H_{1}(f|M) = 1, H_{2}(f|M) = 1,$ and

 $$g = H_{1}(g)\circ f \circ H_{2}(g) |M$$

 for all $g \in U.$

\subsection{Finite determinancy}

We want to see when an isolated singularity of a smooth function is equivalent to that of a polynomial, namely, to a segment of its Taylor series.
It turns out the answer depends on finiteness of an relevant index.

Given an equivalence relation $\sim$ on the set of map-germs $f: (N, S) \rightarrow (P, y)$ (where $N$ and $P$ are manifolds, $S$ is
a finite subset of $N$, and $y$ is a point of $P$), we say $f: (N, S) \rightarrow (P, y)$ is finitely determined if there exists an
integer $k$ such that for any $g: (N, S) \rightarrow (P, y)$ which has the same $k$-jets as $f$ satisfies $f \sim g$.

Let ${\cal F}$ be the set of all $C^{\infty}$ map-germs $f: (N, S) \rightarrow (P, y)$.

Let ${\cal L} = \{$ invertible $C^{\infty}$ map-germs $h^{'}: (P, y) \rightarrow (P, y)\}$
whose group law is composition, and whose action on ${\cal F}$ is given by

$$h^{'} . f = h^{'} \circ f, h^{'} \in {\cal L}, f \in {\cal F}.$$

There is a similar definition of group ${\cal R}$ of invertible $C^{\infty}$ map-germs $h^{'}: (N, S) \rightarrow (N, S)$
with similar action from right on $f \in {\cal F}$.

Let ${\cal A} = {\cal L} \times {\cal R}$ as the product group who acts on $f \in {\cal F}$ from left by ${\cal L}$ and from right by ${\cal R}$.
The orbit of this action is equivalence class of map-germs.

{\bf Map-germ isomorphism}: We say two map-germs $f$ and $g$ are isomorphic if there exist invertible $C^{\infty}$ map-germs
$h: (N, S) \rightarrow (N, S)$ and $h^{'}: (P, y) \rightarrow (P, y)$ such that $h^{'} \circ g \circ h = f$.

We also consider contact equivalence.

Let ${\cal K}$ be the group of invertible $C^{\infty}$ map-germs with composition as group product.
${\cal K}$ acts on ${\cal F}$ with the property $H(\mbox{graph} f) = \mbox{graph} (H \circ f)$.

{\bf Contact equivalence}: The group of germs at $0$ of $C^{\infty}$ diffeomorphisms of $({\bf R}^{n},0)$ acts on the set of $C^{\infty}$ map germs $f:({\bf R}^{n},0)\rightarrow ({\bf R}^{p},0)$ on the right, the group of germs at 0 of diffeomorphisms of $({\bf R}^{n} \times {\bf R}^{p},0)$ of the form $H(x,y)=(x,h_{x}(y))$, with $h_{x}(0)=0$, acts on the graph of $f$ by $H(x, f(x))=(x, h_{x}(f(x)))$.
The equivalence class of this action is Mather's contact equivalence.

Let $N$ be a $C^{\infty}$ manifold and $S$ a subset of $N$. We let $C_{S}$ denote the set of $C^{\infty}$ map-germs
 $(N, S) \rightarrow {\bf R}$. This algebra has a natural ${\bf R}$ structure. Let $m_{S}$ denote the ideal
 in $C_{S}$ consisting of $C^{\infty}$ map-germs $(N, S) \rightarrow ({\bf R}, 0)$. Similar definition applies for $(P, y)$.

 For $f \in {\cal F}$, Mather defined indexes:

 $$d(f, {\cal K}) = \dim_{\bf R} \theta(f) / (t f(B) + f^{*}(m_{y}) \theta(f)),$$

 $$d(f, {\cal A}) = \dim_{\bf R} \theta(f) / (t f(B) + \omega f(A)).$$

  For a map germ $f:({\bf R}^{n},0)\rightarrow({\bf R}^{p},0)$, define $\theta(f)$ as the space of germs of $C^{\infty}$ sections of the bundle $f^{*}T({\bf R}^{p}), tf:\theta({\bf 1}_{{\bf R}^{n}})\rightarrow \theta(f)$ is defined by composing a section of $T({\bf R}^{n})$ with $Tf$,
  $\omega f$ is a vector field along $f$.
  Here $A = \theta({\bf 1}_{(P, y)}), B = \theta({\bf 1}_{(N, S)}), {\bf 1}_{(P, y)}, {\bf 1}_{(N, S)}$ are the identity map-germs of $(P, y)$ and
   $(N, S)$ respectively.

  ${\bf Theorem} \cite{MS3}$: If ${\cal S}$ is any of the groups ${\cal K, A}$, then the necessary and sufficient condition that $f \in {\cal F}$ be
 finitely determined is that $d(f, {\cal S}) < \infty.$

 ${\bf Theorem \cite{MS3}:} f$ is finitely determined with respect to isomorphism if and only if there exists an integer $k$ such that

$$t f(B) + \omega f(A) + m^{l(k)}_{S} \theta(f) \supset m_{S}^{k} \theta(f),$$

for some integer $l(k)$.

\subsection{Algebraic classification of stable germs}

 The problem of classifying stable germs up to isomorphism can be reduced to classifying certain finite dimensional ${\bf R}$-algebra
 up to isomorphism. We define the following algebras:

 $$Q(f) = C_{S}/f^{*}(m_{y}) C_{S},$$

 where $m_{y}$ denotes the unique maximal ideal in $C_{y}$:

 $$Q_{k}(f) = Q(f)/m^{k+1},$$

 where $m$ denotes the intersection of the maximal ideals in $Q(f)$. Let

 $$\bar{Q}(f) = \lim_{\leftarrow} Q_{k}(f).$$

  If $f$ and $f^{'}$ are isomorphic then $Q(f)$ is isomorphic as an ${\bf R}$-algebra to $Q(f^{'})$.
  The converse is also true for stable map-germs.

  Let $x_{1}, ..., x_{s}$ be the distinct points of $S$ and let $x_{1}^{'}, ..., x_{s^{'}}$ be the distinct
  points of $S^{'}$. Let $f_{i} = f|_{(N, x_{i})}: (N, x_{i}) \rightarrow (P, y)$ and define $f^{'}$ similarly.

  ${\bf Theorem} \cite{MS4}$: Suppose that $f$ and $f^{'}$ are stable, that $s = s^{'}, p = p^{'}, n_{i} = n_{i}^{'}$,
  and that there is an isomorphism

  $$Q_{p+1}(f_{i}) \approx Q_{p+1}(f^{'}_{i})$$

  of ${\bf R}$-algebras. Then $f$ is isomorphic to $f^{'}$.


 In summary, stable singularity is uniquely determined by its ring: that is, if the germs are stable and the rings are equivalent,
then they are differentiably equivalent. Thus, the problem of classifying isolated critical points reduces to algebraic problems concerning the action
of finitely dimensional spaces of jets of coordinate maps on jets of map germs.

\subsection{Transversality}

Thom transversality theorem plays an important role in the study of singularity, as seen from the following theorem.

${\bf Theorem} \cite{MS5}$: For a proper $C^{\infty}$ map $f:N \rightarrow P$, assuming $r \ge p+1, k \ge p = \dim P$, the following three conditions are equivalent:

1) $f$ is stable;

2) $f$ is infinitesimally stable;

3) $r j^{k}f$ is transversal to every orbit in $rJ^{k}(N, P)$.

The key notion is multi-jets. Let $N$ and $P$ be manifolds. Let
$N^{r} = \{<x_{1}, ..., x_{r}> \in N^{r}: x_{i} \ne x_{j}, if i \ne j\}$.
Let $\pi_{N}: J^{k}(N, P) \rightarrow N$ denote the projection, where $J^{k}(N, P)$ is the bundle of $k$ jets.
Define $rJ^{k}(N,P) = (\pi_{N}^{r})^{-1}[N^{(r)}],$ where $\pi_{N}^{r}: J^{k}(N, P)^{r} \rightarrow  N^{r}$
is the r-fold Cartisian product of $\pi_{N}$ with itself.
$rJ^{k}(N,P)$ is a fibration over $N^{r} \times P^{r}$. It is the r-fold $k$ jet bundle, or a multi-jet bundle.

If $f: N \rightarrow P$ is a $C^{\infty}$ mapping, $rj^{k}f: N^{(r)} \rightarrow rJ^{k}(N,P)$ is the mapping:
$rj^{k}f(x_{1}, ..., x_{n}) = <j^{k}f(x_{1}), ..., j^{k}f(x_{n})>.$

There is a natural left action of $Diff^{k}N \times Diff^{k}P$ on $rJ^{k}(N,P)$.

By a standard result from consideration of Zariski topology, $rj^{k}f$ is transversal to an orbit in $rJ^{k}(N,P)$
is reduced to the case when the orbit lies in $(\pi_{p}^{r})^{-1} \Delta_{r}$.

Let $f: N \rightarrow P$ be a $C^{\infty}$ mapping, let $x_{1}, ..., x_{r}$ be $r$ distinct points in $N$.
Let $S = \{x_{1}, ..., x_{r}\}$. $\theta(f)_{S}$ is the set of germs at $S$ of $C^{\infty}$ vector fields along $f$.
Let $\theta(N)_{S} = \theta(i_{N})_{S}, tf:\theta(N)_{S} \rightarrow \theta(f)_{S}, tf(\xi) = Tf \circ \xi,$
$\omega f: \theta(P)_{f(S)} \rightarrow \theta(f)_{S}, \omega f(\eta) = \eta \circ f.$

There is a natural identification of ${\bf R}$ vector spaces:

$${\bf T}(rJ^{k}(N,P)_{x})_{z} = \theta(f)_{S}/m^{k+1}_{S} \theta(f)_{S},$$

where $x = (x_{1}, ..., x_{r})$ and $z = rj^{k}f(x).$

$rj^{k}f$ is transversal to $W$ at $x$ if, and only if,

$$t f(\theta(N)_{S}) + \omega f(\theta(P)_{f(S)}) + m_{S}^{k+1} \theta(f)_{S} = \theta(f)_{S}.$$

 \subsection{Density of topological stable maps}

How generic are stable maps in the space of maps? Thom conjectured that if one is allowed all homeomorphisms as coordinate changes, stable maps are dense
in the space of maps. In \cite{MS7} a rigorous proof is provided. More importantly it establishes the foundation for stratification
of singular spaces.

  {\bf Theorem:} Let $X$ be a complex subvariety of a smooth variety $M$, then $M$ admits a Whitney stratification such that
  $X$ is a finite union of strata.

  A Whitney stratification satisfies Whitney regular condition.

  {\bf Whitney regular condition}: Let ${\cal S} = \{S_{i}| i \in I\}$ be a set of stratas. Let $X, Y \in {\cal S}, X \cap Y = \phi$.
  Let $x \in \bar{Y} \cap X, (Y, X, x)$ is Whitney regular if, for any sequence
  $(x_{n}, y_{n}) \in X \times Y, x_{n} \rightarrow x, y_{n} \rightarrow x, v_{n} = \frac{1}{|y_{n} - x_{n}|}(y_{n} - x_{n}) \rightarrow v$
  $T_{y_{n}} \rightarrow T$, then $v \in T$.

 \subsection{Classification of holomorphic germs: Mather-Yau's theorem}

 Mather's works above establish the foundation of singularity theory. The techniques developed for singularities of smooth germs
 can also be used for singularities of complex germs.

 Let ${\cal O}_{n+1}$ denote the ${\bf C}$-algebra of germs at the origin of holomorphic functions $f: ({\bf C}^{n+1},0) \rightarrow {\bf C}$,

 If $V:=\{f=0\}$ is the hypersurface germ at the origin defined by some $f \in {\cal O}_{n+1}, f(0)=0$,

 Associate to $V$ two ${\bf C}$-algebras:

 $A(V)={\cal O}_{n+1}/(f, \frac{\partial f}{\partial z_{0}}, ...,\frac{\partial f}{\partial z_{n}}) {\cal O}_{n+1},$

 $B(V)={\cal O}_{n+1}/(f, z_{i} \frac{\partial f}{\partial z_{j}}){\cal O}_{n+1}$ with $0 \le i,j \le n$.

  Two hypersurface germs $(V,0)$ and $(W,0)$ in ${\bf C}^{n+1}$ with $V / 0$ nonsingular are biholomorphically equivalent if and only if $A(V)$ and $A(W)$ (resp. $B(V)$ and $B(W)$) are isomorphic as ${\bf C}$-algebras.

 Here $A(V)$ is the moduli algebra of $V$ (the germ at $0$ of $A(V)$ is the base space for the miniversal deformation of the hypersurface $V$).

  Let $J^{k}$ denote the ${\bf C}$-vector space of $k$-jets at the origin of elements of $O_{n+1}$. Let $K^{k}$ denote the Lie group of $k$-jets
  at the origin of members of $K$. Since $K$ acts on $O_{n+1}$, we have that $K^{k}$ acts on $J^{k}$. For $f \in O_{n+1}$, let $f^{(k)} \in J^{k}$
  denote the $k$-jet of $f$ at the origin. We say $f$ is $k$-determined with respect to $K$ if $g \in O_{n+1}$ and $g^{k} \in K^{k}f^{k}$ imply
  $g \in Kf$. We say $f$ is finitely determined with respect to $K$ if it is $k$-determined with respect to $K$ for some positive integer $k$.

  The proofs depends on techniques developed by Mather in singularity theory: group actions on jet spaces, finite determinacy and Mather's
  sufficient conditions for a connected submanifolds to be contained in an orbit.

\subsection{Holomorphic map-germs, continued}

Arnold also did great works in singularity theory, especially in complex singularities. He also use singularity theory
to explain experiments such as caustics which is the original goal of Thom.

 As above classification of singularities reduced to a finite dimensional problem with a finite dimensional 
 Lie group (coming from finite segment of jets of isomorphism of map germs) acting on finite segment of jets.
 
 Let $G$ be a Lie group acting on a variety $X$, the modality $m$ of a point $x \in X$ is the smallest number
such that a sufficiently small neighborhood of $x$ is covered by a finitely many $m$-parameters of orbits.
 The modality of the germ of a function at a critical point with critical value $0$ is defined as the modality
of its sufficient jets in the space of jets of functions with critical point $O$ and critical value $0$.

 Two germs are stably equivalent if they become equivalent after direct addition of non-degenerate quadratic forms.
 Germs can be classified through its associated algebra such as $A(V)$. Adding a non-degenerate quadratic form
 implies adding a trivial factor of the associated algebra.

 Up to stable equivalence, the germs with $m=0$ are exhausted by the following list:

 $A_{k}: f(x) = x^{k+1} , k \ge 1;$

 $ D_{k}: f(x, y) = x^{3}y + y^{k-1}, k \ge 4;$

 $E_{6}: f(x,y) = x^{3} + y^{4};$

 $ E_{7}: f(x, y) = x^{3} + x y^{3};$

 $ E_{8}: f(x,y) = x^{3} + y^{5}.$

 They are called {\bf ADE singularities} which are related to Dynkin diagrams.

 With each isolated critical point of a holomorphic function, consider a sufficiently small ball with centre at the critical point.
 The part of the level set inside the ball and sufficiently close to the critical set is a smooth manifold $V$ with boundary $\partial V$.
 According to Milnor \cite{Milnor1}: $\partial V$ is homotopic equivalent to $\mu$ middle dimensional vanishing cycles, $S^{n-1} \bigvee ... \bigvee S^{n-1}$.

 Let $\mu = \dim M(f):= C[z_{1}, ..., z_{n}] / \{\frac{\partial f}{\partial x_{1}}, ..., \frac{\partial f}{\partial x_{1}}\}$.
 Then $H_{n-1}({\partial V}, {\bf Z}) = {\bf Z}^{\mu}$.

 The intersection index defines the quadratic form of the singularity.
 The quadratic form of simple singularities $A, D, E$ are given by the Dynkin diagrams.

 There is a miniversal deformation of the function $\lambda \rightarrow f + \lambda_{1} e_{1} + ... + \lambda_{\mu} e_{\mu}$.
 Every deformation of $f$ is equivalent to the one induced from a versal deformation. The level bifurcation set is formed by $\lambda$
 for which $0$ is the critical value. The complement to the level bifurcation set is the base of a fibration whose fiber is a non-singular
 level manifold of $f$. The action of the fundamental group of the base on the homology of a fibre is monodromy of the singularity.

 {\bf Remark:} In the case of Yang-Yang functions the monodromy can be used to construct knot invariants. It reproduces Jones polynomials.

\section{Mather's work on Hamiltonian Dynamical Systems}

  Since late seventies John devoted most of his time to the development of variational methods for Hamiltonian systems.
  He started his works on dynamics of two dimensional area preserving mappings and later he extends it to many degrees of freedom.
  As an application it solves some of the problems on Arnold diffusion.

\subsection{Aubry-Mather set}

 Let $f: (x, y) \rightarrow (x^{'}, y^{'}), x \in {\bf R}/{\bf Z}, y \in {\bf R}$ be an area-preserving diffeomorphism.
 Now $y^{'} dx^{'} - y dx$ is closed. Assuming it is exact then it equals to $dh(x, x^{'}), h$ is a generating function.
 We have $y = - \partial_{1}h(x, x^{'}), y^{'} = \partial_{2}h(x, x^{'})$.
 We further assume it obeys the twist condition: $\frac{\partial y}{\partial x^{'}} > 0.$
 Then the map $(x, y) \rightarrow (x, x^{'})$ is a coordinate change.

 Let us consider configuration $(..., x_{1}, x_{2}, ..., x_{n}, ...) \in ({\bf R}/{\bf Z})^{\infty}$.
 Over each finite segment $(x_{i_{1}}, x_{i_{1}+1}, ..., x_{i_{2}})$ with fixed end points $x_{i_{1}}, x_{i_{2}}$
 one can consider a variational principle:
 $$H(x_{i_{1}}, x_{i_{1}+1}, ..., x_{i_{2}}) = \Sigma_{i=i_{1}}^{i_{2}} h(x_{i}, x_{i+1}).$$
 It is easy to see that critical point of $H$ gives an orbit
 $$(x_{i_{1}}, y_{i_{1}}), ..., (x_{i_{2}}, y_{i_{2}})$$

 where $y_{i} = - \partial_{1} h(x_{i}, x_{i+1}) = \partial_{2} h(x_{i-1}, x_{i}).$

  By minimizing $\Sigma_{i_{1}}^{i_{2}} h(x_{i}, x_{i+1})$ for each finite segment it gives a orbits. Such an orbit has a well-defined rotation number.
 For an area-preserving twist map and a rotation number, the totality of such minimal orbits is called the Aubry-Mather set \cite{M-HD0}.

 As an example one may consider minimal geodesics over a surface. Orbits of Aubry-Mather set correspond to minimal
geodesics. Such minimal geodesics was constructed by Hedlund \cite{H}. One can assign such minimal geodesics with various rotation numbers
such as $p/q, p/q^{+}, p/q^{-}$. There are periodic orbits of hyperbolic type and homoclinic orbits of those orbits.

 One can define Peierl's energy barrier:

 $P_{\omega}(a) = \min \{\Sigma_{i=-\infty}^{\infty} h(x_{i}, x_{i+1}) - \Sigma_{i=0}^{\infty} h(x_{i}^{+}, x_{i+1}^{+})
- \Sigma_{i=-\infty}^{0} h(x_{i}^{-}, x_{i+1}^{-}): x_{0} = a\}$

 with $(x_{0}^{+}, x_{1}^{+}, ..., ), (..., x_{-1}^{-}, x_{0}^{-})$ given action minimizing configurations of rotational number $\omega$.

 We have $P_{\omega}(a) = 0$ if and only if there is a minimizing orbit through $a$.
 Let $M_{\omega}$ be the union of all minimizing orbits, it is Aubry-Mather set.

 There exists an invariant circle with rotation number $\omega$ if and only if $P_{\omega}(a) = 0$ for every $a$.

\subsection{Mather's connecting orbits in Birkhoff region of instability}

 Existence of invariant circle of a rational rotation number is not generic. The set of rotation numbers such that there are invariant circles
 with such rotational numbers is closed. Its complement consists of countably many intervals. For each such an interval $I = (a, b)$, to the end points
 $a, b$, there are invariant circles $\Gamma_{a}, \Gamma_{b}$. There is no invariant circle with rotational number $c \in (a, b)$.
 The phase space $R$ bounded by $\Gamma_{a}, \Gamma_{b}$ is called a {\bf Birkhoff region of instability}.

 Question 1: Are the rotation numbers of boundary invariant curves always irrational?

 Question 2: Are there self-similar property for the boundary invariant curves?

 {\bf Mather's connecting orbits in Birkhoff region of instability \cite{M-HD2}:} Let $R$ be a Birkhoff region of instability, bounded by invariant curves $\Gamma_{a}, \Gamma_{b}$. Suppose $a < \alpha, \omega < b$, then there is an orbit whose $\alpha$-limit set is $M_{\alpha}$ and whose $\omega$-limit set is $M_{\omega}$.

\subsection{Minimal measures for several degree of freedoms}

 Mather extends the above works to several degree of freedoms. A crucial concept is that of invariant measures.
 Here is he setting for several degrees of freedoms.

 Let $L: TM \times {\bf T} \rightarrow {\bf R}$ be a Lagrangian, here ${\bf T} = {\bf R} / {\bf Z}$. Assume it has the following properties:

 1) Positive definiteness: $L_{\dot{x}, \dot{x}}$ is positive definite.

 2) Superlinear growth: Let $||.||$ denote a Riemannian metric on $M$. Then

 $$L(v, \theta) /||v|| \rightarrow + \infty, \mbox{ as } ||v|| \rightarrow \infty,$$

 where $v$ ranges over $TM$ and $\theta \in {\bf T}$.

 3) Completeness of the Euler-Lagrangian flow: Every maximal trajectories of the Euler-Langrange flow is defined for all time.

  Our goal is to understand the space of orbits of the corresponding Euler-Langrange flow.
  It is better to look at the space of invariant probability measures on the phase space.

 Let $\mu$ be a probability measure on the phase space $TM \times {\bf T}$, invariant under the Euler-Lagrange flow,
 the space of invariant probability measures is a convex set.
 The average action of $\mu$ is $A_{L}(\mu) = \int L d \mu$, it is an affine function of $\mu$.

 For each invariant probability measure we can associate a rotation vector $\rho(\mu) \in H_{1}(M, {\bf R})$.
 For $c \in H^{1}(M, {\bf R})$, let $\lambda_{c}$ be a closed 1-form on $M$ representing $c$.
 The rotation vector $\rho(\mu)$ is uniquely determined by $<c, \rho(\mu)> = \int \lambda_{c} d \mu$.

{\bf Minimal measures}

 Consider $U_{L} = \{(\rho(\mu), A(\mu)): \mbox{ $\mu$ is an invariant probability measure with $A(\mu) < \infty$} \}$.
 $U_{L}$ is a convex set and $U_{L} \subset H_{1}(M, {\bf R}) \times [B, \infty)$, therefore $U_{L}$ is the epigraph of a convex function
 $\beta = \beta_{L}: H_{1}(M, {\bf R}) \rightarrow {\bf R}$.
 For $h \in H_{1}(M, {\bf R}), \beta(h)$ is the minimal average action of the rotation vector $h$.

 Let $\alpha: H^{1}(M, {\bf R}) \rightarrow {\bf R}$ be the conjugate of $\beta$, i.e. for $c \in H^{1}(M, {\bf R})$,

 $$\alpha(c) = \max \{<c, h> - \beta(h): h \in H_{1}(M, {\bf R})\}.$$

 Consider $A_{c}(\mu) = A(\mu) - <c, \rho(\mu)> = \int (L - \lambda_{c}) d \mu$.
 It is easy to see that $L$ and $L - \lambda_{c}$, with $\lambda_{c}$ a closed one-form, have same extremals.

 {\bf c-minimal}: It is an invariant probability measure to minimize $\int (L - \lambda_{c}) d \mu$
over the set of invariant probability measures. We denote ${\cal M}_{c}$ the set of all c-minimal measures.
Let $M_{c}$ be the closure of the union of the supports of $\mu \in {\cal M}_{c}$.

 We now extend {\bf variational principle} to several degree of freedom.

Let $\tilde{M}$ be the covering space of $M$ such that $\pi_{1}(\tilde{M}) = ker ({\cal H}: \pi_{1}(M) \rightarrow H_{1}(M, {\bf R}))\}$

Define $h: \tilde{M} \times \tilde{M} \rightarrow R$ as $h(\tilde{m}, \tilde{m}^{'}) = \min \int_{0}^{1} L(d \gamma(t), t) dt$,
where the minimum is taken over all curves $\tilde{\gamma}: [0,1] \rightarrow \tilde{M}$ such that $\tilde{\gamma}(0) = \tilde{m}, \tilde{\gamma}(1) = \tilde{m}^{'}$.

 A configuration is a bi-infinite-sequence $(..., \tilde{m}_{i}, ...), \tilde{m}_{i} \in \tilde{M}$.
A segment is a finite sequence $(\tilde{m}_{a}, ..., \tilde{m}_{b})$. We define a function over such a segment,
 $h(\tilde{m}_{a}, ..., \tilde{m}_{b}) = \Sigma_{i = a}^{b-1} h(\tilde{m_{i}}, \tilde{m}^{'}_{i+1}).$

 A segment is minimal if $h(\tilde{m}_{a}, ..., \tilde{m}_{b}) \le h(\tilde{m}_{a}^{'}, ..., \tilde{m}_{b}^{'})$,
for any segment with two end points fixed. A configuration is minimal if every segment of it is minimal, it also corresponds
to a minimal orbit in $\tilde{M}$.
For a $\tilde{M}$-minimizer $\gamma$, there exists a $c \in H^{1}(M, {\bf R})$ such that every limit measure of $\gamma$ minimizes $A_{c}$.

{\bf Generalization of Peierls barriers}

 For $m, m^{'} \in M$, set $h_{c}(m, m^{'}) = \min \int_{0}^{1}(L - \lambda_{c}) (d \gamma(t), t) dt - \alpha(c)$.
 A segment $(m_{a}, ..., m_{b})$ of an $M$-configuration is $c$-minimal if
 $h_{c}(m_{a}, ..., m_{b}) \le h_{c}(m_{c}^{'}, ..., m_{d}^{'})$, for any other segment with the same ends.
 The n-fold conjunction of $h_{c}$ is $h_{c}^{n}(\xi, \eta)  = \min \{ \Sigma_{i=0}^{n-1} h_{c}(m_{i}, m_{i+1}): m_{0} = \xi, m_{n} = \eta \}$.

 Set $h_{c}^{\infty}(\xi, \eta) = \lim \inf_{n \rightarrow \infty} h_{c}^{n}(\xi, \eta), \xi, \eta \in M.$
 Set $B_{c}(\xi) = h_{c}^{\infty}(\xi, \xi)$, it is a generalization of Peierls barrier.
 The zero set of $B_{c}(\xi)$ is the support of $c$-minimals.

{\bf Generalization of Birkhoff region of instability to several degree of freedoms}

  Set $B_{c}^{*}(m) = \min \{h_{c}^{\infty}(\xi, m) + h_{c}^{\infty}(m, \eta) - h_{c}^{\infty}(\xi, \eta): \xi, \eta \in \Sigma_{c}^{0}\}$.
 Let $W_{L} = \{c \in H^{1}(M, {\bf R}):$ there exists an open neighborhood $U$ of $\{B_{c}^{*} = 0\}$
 in $M$ such that the inclusion map  $H_{1}(U, {\bf R}) \rightarrow H_{1}(M, {\bf R})$ is the zero map$\}$.
 It is a generalization of Birkhoff region of instability to several degree of freedoms.

 {\bf Mather's connecting orbits in Birkhoff region of instability}:
 Suppose $c_{0}$ and $c_{1}$ are in the same connected component of $W_{L}$. Then there is a trajectory of the Euler-Lagrangian flow
 whose $\alpha$-set lies in $M_{c_{0}}$ and whose $\omega$-limit set lies in $M_{c_{1}}$.
 
  In fact connecting orbits can be constructed for more general case. Given $c \in H^{1}(M, {\bf R})$, we define
  
  $$V_{c} = \cap_{U} \{ i_{U*} H_{1}(U, {\bf R}) | U \mbox { is a neighborhood of } \{B_{c}^{*} = 0\} \}.$$
  
  Here $i_{U}: U \rightarrow M$ is the inclusion map. Define
  
  $$V_{c}^{\bot} = \cup_{U} \{ker i_{U}^{*}: U \mbox{ is a neighborhood of $\{B_{c}^{*} = 0\}$}\}.$$
  
  We say that $c_{0}, c_{1} \in H^{1}(M, {\bf R})$ are $c$-equivalent if there exists a continuous curve 
  $\Gamma: [0, 1] \rightarrow M$ such that $\Gamma(0) = c_{0}, \Gamma_{1} = c_{1}$, and for each $t_{0} \in (0, 0)$,
  there exists $\delta > 0$ such that $\Gamma(t) - \Gamma(0) \in V_{\Gamma(t_{0})}^{\bot}$ whenever $t \in [0, 1]$
  and $|t - t_{0}| < \delta. M_{c_{0}}, M_{c_{1}}$ can be connected by an orbit whenever they are $c$-equivalent.
  We may call this approach Mather's mechanism in connecting different regions in the phase space.
  
  His proof relies on {\bf The method of changing Lagrangians}:

 Choose a sequence ..., $S_{i}$, ... of hypersurfaces of $M \times {\bf T}$;
a sequence ..., $T_{i}$, ... of positive numbers; a sequence ..., $\eta_{i}$, ... of closed one-forms.
 For $(\sigma, \tau) \in S_{i}, (\sigma^{'}, \tau^{'}) \in S_{i+1}$, set

 $h_{i}((\sigma, \tau), (\sigma^{'}, \tau^{'})) = \inf \{ \int_{a}^{b} (L - \eta_{i})(\theta, \dot{\theta}, t) dt \}$
 where the infinum is taken over all curves $\theta: [a, b] \rightarrow M$ such that
 $a = \tau_{i} \mod 1, b = \tau_{i+1}, \theta(a) = \sigma, \theta(b) = \sigma^{'}, b-a \ge T_{i}.$
 
 To join several minimal segments one needs Peierls energy barrier at the joint region so that
 minimizers do not hit the constraints.

 \subsection{Applications to Arnold diffusion}

 Consider an near integrable system with the following Lagrangian:

 $L(\theta, \dot{\theta}, t) = l_{0}(\dot{\theta}) + \epsilon P(\theta, \dot{\theta}, t), \theta = (\theta_{1}, ..., \theta_{n}) \in {\bf T}^{n}, {\bf T} = {\bf R}/{\bf Z}.$

 When $\epsilon = 0$ the system is integrable. Solutions are periodic or quasi-periodic orbits
 of frequencies $\partial l_{0} / \partial \dot{\theta}.$ The full phase space is foliated by invariant tori.

 When $\epsilon$ is small and $\partial^{2} l_{0}$ is non-degenerate invariant tori with Diophantine frequency vectors
 survived after Kolmogorov-Arnold-Moser \cite{A2}. When $n=1$ those tori divided the phase space into different regions
 and motions are confined by those tori. When $n = 2$ Arnold constructed an example with motions drifted away
 significant amount in action. The Hamiltonian of his example reads:

 $$H =  \frac{p_{1}^{2}}{2} + \frac{p_{2}^{2}}{2} + \epsilon(\cos q_{1} - 1)(1 + \mu (\cos q_{2} + \cos t)) .$$

 If $\mu = 0$ it is a decoupled rotator-pendulum system. Motions of rotator are simply rotations.
 Pendulum has a hyperbolic fixed point with stable and unstable manifolds. Motions of rotators
 are normally hyperbolic invariant cylinders. Under perturbation those normally hyperbolic cylinders
 still survive. In the normal direction one can drift significant amount in action.

 People asked if Arnold's phenomenon is generic for near integrable systems. In Lagrangian formalism Mather gave a
 definition of Arnold diffusion. Set
 $ osc_{(\theta, J)} \dot{\theta} = \sup \{||\dot{\theta}(t_{0}) - \dot{\theta}(t_{1})||: t_{0}, t_{1} \in J\}$
associated to a trajectory $\theta: J \rightarrow {\bf T}^{n}.$

 {\bf Arnold diffusion:} For fixed $l_{0}$ and for some $\delta > 0$, there exists a set of typical small perturbations $P$
such that for any perturbation in the set there exists a trajectory $(\theta, J)$ for which

 $ osc_{(\theta, J)} \dot{\theta} > \delta.$

 In Hamiltonian formalism general near integrable system reads: $H(p, q) = h (p) + \epsilon P(p, q)$.

 A frequency vector $\omega(p) = \partial h /\partial p $ is called a resonance vector if
 there is a vector of integers $\vec{k}$ such that $\{<\vec{k}, \omega(p)> = 0\} \cap \{h(p) = 0\}$.
 With a path of resonance vector one can kill a coordinate in the phase space and put normal form of Hamiltonian
 with one fewer coordinate. This way one is at a similar
 position as Arnold's example. It is an integable system coupled with a twist map. Normally hyperbolic invariant cylinders
 persists after a small perturbation. One can apply Mather's variational methods to construct connecting
 orbits in different regions in the normal direction. For a path of single resonance without double resonance,
 or a priori unstable systems, several authors succeeded in constructing such connecting orbits \cite{Xia, Cheng1, BKZ}.

However in general in the path of a single resonance one often meets with double resonance frequency vectors, i.e.
$\{<\vec{k} , \omega(p)> = <\vec{k}^{'}, \omega(p)> = 0\} \cap \{h(p) = 0\}$. In this case the normal form reads:

$$G(X, Y) = \frac{1}{2} <Ay, y> + V(x), (x, y) \in T^{*} {\bf T}^{2}, \max V = 0.$$

In the case of double resonance one needs to do more work. One can apply variational methods to
construct normally hyperbolic invariant cylinders and connecting orbits \cite{Cheng2}, see also \cite{KZ, Marco}. 
For a nice introduction to the current stage of Arnold diffusion please see \cite{Cheng3} in this volume. 
Finally one is able to prove:

{\bf Mather's Last Theorem:} (\cite{M-HD3, M-DS4, M-DS5}) Set $n = 2$. Let $\Omega_{1}, ..., \Omega_{k}$ be open subsets of $B^{n}$.
 There exists a non-negative continuous function $\delta$ defined on ${\cal P}^{3} \times {\cal L}^{3}$ such that:
 For any $l_{0} \in {\cal L}^{3}, \{P \in {\cal P}^{3}: \delta(P, l_{0}) > 0\}$ is dense in ${\cal P}^{3}$,
 There exists an open, dense subset $W$ of $\{(\epsilon, l_{0}, P) \in {\bf R}_{+} \times {\cal L}^{3} \times {\cal P}^{3}:
 0 < \epsilon < \delta(P, l_{0}) \}$ such that for $(\epsilon, l_{0}, P) \in W$, there exists a trajectory
of $L = l_{0} + \epsilon P$ that visits the $\Omega_{i}$'s in any pre-assigned order. 
Here ${\cal P}^{r}$ denote the topological space of $C^{r}$ functions $P: {\bf T}^{n} \times B^{n} \times {\bf T} 
\rightarrow {\bf R}$ such that $||P||_{3} = 1,$ provided with the $C^{r}$-topology. 
${\cal L}$ denote the topological space of $C^{r}$ functions $l_{0}: B^{n} \rightarrow {\bf R}$
such that $d^{2} l_{0} > 0$, provided with the $C^{r}$ topology.

\section{Mather-Chern classes over complex varieties}

    Characteristic classes were constructed by a number of great mathematicians E. Stieffel, H. Whitney, L. Pontrajin and S. S. Chern \cite{Milnor}.
The aim is to provide an invariant to distinguish bundles over some spaces. When the space is smooth it has a tangent bundle.
When the space is singular tangent space may not exist. That is the problem to construct Chern classes for a singular
complex variety.

    In 1969, P. Deligne and A. Grothendick conjectured existence and uniqueness of Chern classes for singular algebraic varieties.
The conjecture was proved by B. MacPherson in 1972. One of the fundamental ingredients of the MacPherson construction is the
Mather class, defined using the Nash transformation. The MacPherson class is a combination of Mather classes, the coefficients
being defined using the local Euler obstruction.

  In 1965, two constructions of Chern classes for singular varieties have been published. The one was published by Wu Wen Tsun,
  in Chinese, the other was published by M. Schwartz, in French.

    The relation between complex Wu classes and Schwartz classes appeared thirty years later, through the works of B. MacPherson
    and Jianyi Zhou. In short, Wu class is equivalent to Mather class and Schwartz class is equivalent to MacPherson class.

    \subsection{Mather classes}

    Let us consider the Nash transformation $\mu: \hat{X} \rightarrow X$ and the Nash bundle $\hat{T}$ with basis $\hat{X}$,
    that is the bundle where fibre at $(x, P) \in \hat{X}$ is

    $$\{(x, \nu, P)| \nu \in P\}.$$

    The Mather classes of $X$ are defined by

    $$C^{M}(X) = \nu_{*}(c^{*}(\hat{T}) \cap [\hat{X}]),$$

    where $c^{*}(\hat{T})$ is the total cohomology Chern class of the bundle $\hat{T}$ and the Poincare homomorphism
    $H^{*}(\hat{X}) \rightarrow H_{d-*}(\hat{X})$, cap-product by the fundamental class $[\hat{X}]$ is not necessarily an isomorphism.

    R. Piene \cite{Piene} provided an expression of Mather classes in terms of polar varieties. The $k$-th polar
    variety of $X$ relative to a linear subspace $L_{k}$ of codimension $d-k+2$ in ${\bf C P}^{n}$ is defined as

    $$M_{k} = \mbox{closure of} \{x \in X_{reg}| \dim(T_{x}(X_{reg}) \cap L_{k}) \ge k-1\}.$$

    For a linear subspace $L_{k}$ in general position, $M_{k}$ represents a class of rational equivalence of codimension $k$ in $X$,
    denoted by $[M_{k}]$.

    Let us denote by ${\cal L}$ the restriction of the hyperplane bundle to $X$,

    $${\cal L} = {\cal O}_{{\bf C P}^{n}}(1)|_{X},$$

    then one has

    $$C^{M}_{d-k}(X) = \Sigma_{i=0}^{k} (-1)^{i} \binom{d-i+1}{d-k+1} c_{1}({\cal L}^{k-i}) \cap [M_{i}].$$

    \subsection{Mather classes and Wu classes are equivalent}

    Let $X$ be a $d$-dimensional algebraic complex projective variety in ${\bf C P}^{n}$.

    {\bf Definition:} One says that $\xi$ is a generic point of $X$ if, for an extension $\tilde{\bf C}$ of ${\bf C}$,
    then $\xi$ is a point of $\tilde{\bf C}{\bf P}^{n}$ such that $\tilde{X} \cap {\bf C P}^{n} = X$, where
    $\tilde{X} = \bar{\{\xi\}}$ is the Zariski closure of $\{\xi\}$ in $\tilde{\bf C}{\bf P}^{n}$.

    A proper speciafication of $\xi$ is a solution, on ${\bf C}$, of equations satisfied by $\xi$.

    A generic point $\xi$ of $X$ is a simple point, that means that if $\{f_{i}\}$ is a finite number of
    polynomials defining $X$, then the Jacobian matrix $((\frac{\partial f_{i}}{\partial x_{j}})(\xi))$
    has maximal rank $n-d$. The projection tangent space $T_{\xi}(\tilde{X})$ of $\tilde{X}$ at $\xi$ is
    well defined in $\tilde{\bf C}{\bf P}^{n}$.

    {\bf Definition:} The variety $\hat{X}_{W}$ in $G$ whose $(\xi, T_{\xi}(\tilde{X})$ is a generic point
    is called Wu transformation of $X$.
    Here $G = \{(x, P)| x \in P, P \mbox{ is a d-plane in } {\bf C P}^{n}\}$ is an algebraic manifold of dimension
    $m = n + d(n-d).$

    The restriction to $\hat{X}_{W}$ of the projection of $G$ on ${\bf C P}^{n}$ is a birational map $\nu: \hat{X}_{W} \rightarrow X$.

    Let $0 \le b_{0} < b_{1} < ... < b_{d} \le n$ be a set of integers and $a$ is one of $b_{i}$. Denote:

    $[a/b_{0}, b_{1},..., b_{d}]$ is a cycle in $G$ with dimension $a + \Sigma^{'}_{i}(b_{i} - i)$,\
    where $\Sigma_{i}^{'}$ is the sum on all $i$ such that $b_{i} \ne 0$.

    Using Wu transformation Wu defined rational equivalent classes of $W_{s}([a/b_{0}, b_{1}, ..., b_{d}])$ valued in
    Chow ring of $X, A_{*}(X)$. The Wu class is defined by

    $$C^{W}_{d-s}(X) = \Sigma_{i=0}^{s} (-1)^{i} \binom{d-i+1}{d-s+1} W_{s}([s-i/0, ..., d-i, d-i+2, ..., d+1]).$$

    {\bf Theorem:} The Wu and Nash transformation coincide: $\hat{X}_{W} = \hat{X}_{N}.$

    As a consequence Mather classes and Wu classes coincide.



\section{Foliations: Mather-Thurston Theorem}

Haefliger \cite{Ha} has shown that the cohomology of $B\Gamma$ measures the obstruction
to finding a foliation on an open manifold with a given normal bundle.

In \cite{M-Fo1, M-Fo2} it is proved that $H_{i}(BG)=H_{i}(\Omega B\Gamma)$.
Naturally it has deep applications to theory of foliations.

Thurston \cite{Th1} has found a marvelous generalization of the main result to
higher codimensions.

 \section{Mather's other works}

 John has also made deep contributions to celestial mechanics and characterization of prime ends with applications to conformal mapping over the boundary.
 We will not make any further introductions to them.

\section{Summary}

 John Mather is a great scholar who was dedicated to mathematics in his whole life.

 His works in mathematics can be characterized as original and foundational. He laid out the foundation of singularity theory while he was a graduate student. He also laid out the foundation of modern Hamiltonian dynamical systems. Those fields became main stream in mathematics and it attracts many talents to pursue. His other works on characteristic classes, foliations, celestial mechanics, prime ends of conformal mappings are of the same quality with great influence in mathematics.

\end{document}